\documentclass[12pt]{article}
\usepackage{epsf,epic,eepic,float,fullpage}
\usepackage{graphicx}
\usepackage{float}
\restylefloat{figure}

\def\RR{\hbox{I\kern-.2em\hbox{R}}}
\newcommand{\qed}{\hbox to 0pt{}\hfill$\rlap{$\sqcap$}\sqcup$
\vspace{3mm}}

\title{Stability  of the Second Order Delay Differential Equations 
with a Damping Term}
\author{Leonid Berezansky $^{1}$ \\Department of Mathematics,
Ben-Gurion University of the Negev, \\
Beer-Sheva 84105, Israel,
\\
Elena Braverman $^{2}$ \\
Department of Mathematics and Statistics, University of Calgary, \\
2500 University Drive N.W., Calgary, AB T2N 1N4, Canada\\
and   Alexander Domoshnitsky  $^{3}$ \\
Department of Mathematics and Computer Science,\\
The College of Judea and Samaria,
Ariel 44837, Israel}
\date{}
\begin{document}
\maketitle

\begin{abstract}
For the delay differential equations 
$$
\ddot{x}(t) +a(t)\dot{x}(g(t))+b(t)x(h(t))=0,~ g(t)\leq t,~h(t)\leq t,
$$
and
$$
\ddot{x}(t) +a(t)\dot{x}(t)+b(t)x(t)+a_1(t)\dot{x}(g(t))+b_1(t)x(h(t))=0
$$
explicit exponential stability conditions are obtained.
\end{abstract}


{\bf AMS(MOS) subject classification:} 34K20.

\thispagestyle{empty}

\footnotetext[1]{Partially supported by Israeli Ministry of Absorption
and the Israel Science Foundation (grant No. 828/07)}
\footnotetext[2]{Partially supported by the NSERC Research Grant}
\footnotetext[3]{Partially supported by the Israel Science Foundation 
(grant No. 828/07)}

\section{Introduction}

This paper deals with the scalar linear delay 
differential equation of the second order with a damping term
\begin{equation}
\label{1}
\ddot{x}(t) +a(t)\dot{x}(g(t))+b(t)x(h(t))=0
\end{equation}
and the equation which also involves nondelay terms
\begin{equation}
\label{2}
\ddot{x}(t) +a(t)\dot{x}(t)+b(t)x(t)+a_1(t)\dot{x}(g(t))+b_1(t)x(h(t))=0.
\end{equation}
Such linear and nonlinear equations attract attention 
of many mathematicians due to their 
significance in applications.
We mention here the monographs of 
Myshkis~\cite{M},
Norkin~\cite{N},  Ladde, Lakshmikantham and  Zhang~\cite{LLZ}, 
Gy\"{o}ri and  Ladas~\cite{GL},  Erbe,  Kong and Zhang~\cite{EKZ}, 
Burton~\cite{Bur}, Kolmanovsky and
Nosov~\cite{KN}  and references therein.

In particular, Minorsky~\cite{Min} in 1962 considered the problem of 
stabilizing the rolling of a
ship by the "activated tanks method" in which ballast water is pumped from one position to
another. To solve this problem he constructed several delay differential
equations with damping of the form (\ref{1}) and (\ref{2}).

In spite of  obvious importance in applications, there are only few 
papers on delay differential equations with damping. 

In \cite{CS} the authors considered  autonomous equation (\ref{2}) 
and obtained stability results using analysis of the roots of
the characteristic equation.

In \cite{Bur} stability of the autonomous equation
\begin{equation}
\label{3}
\ddot{x}(t) +a\dot{x}(t)+bx(t-\tau)=0
\end{equation}
was studied using Lyapunov functions. It was demonstrated that
if $a>0, b>0$ and $b\tau<a$, then equation (\ref{3}) is exponentially 
stable. 
Other results obtained by Lyapunov functions method can be found in 
\cite{BH,Z}.

In  \cite{BurFur} Burton and Furumochi applied  fixed point theorems to equation
(\ref{1}) and obtained new stability results.
In particular, the equation
$$
\ddot{x}(t) +\frac{1}{3}\dot{x}(t)+\frac{1}{48}x(t-16)=0 
$$  
is exponentially stable, where the condition $b\tau<a$  
does not hold. Here $b\tau = a$.

In \cite{Bur2,Bur3}  some other stability conditions were obtained 
by the fixed point method
for (\ref{1}) in the case $g(t)\equiv t$, $h(t)=t-\tau$.

To the best of our knowledge, there is only one paper \cite{BD} where 
stability of  the 
general nonautonomous equation  (\ref{1}) was investigated.
In \cite{BD}  the authors  applied the W-method
\cite{AS} which is based on the application of the Bohl-Perron type 
theorem (Lemma 3 of the present paper).

Note also the paper \cite{D} where nonoscillation of systems 
of delay differential equations was considered and on this basis 
several results on nonoscillation and exponential stability of second 
order delay differential equations were obtained.

Here we will employ the method of \cite{BD} and consider equation 
(\ref{2}), which was not studied in \cite{BD}. Furthermore, we will use a 
new approach (also based on a Bohl-Perron type theorem) and obtain 
sharper stability results for equation (\ref{1}) than in \cite{BD}. 

In particular, for (\ref{3}) we  obtain the same stability 
condition $b\tau<a$ as in \cite{Bur}, but our results are applicable to 
more general nonautonomous equations as well.

\section{Preliminaries}

We consider the scalar second order delay differential equation (\ref{1})
under the following conditions:

(a1) $a(t), b(t),$ are Lebesgue measurable and  
essentially bounded functions on $[0,\infty)$;

(a2) $g:[0,\infty)\rightarrow \RR, h:[0,\infty)\rightarrow \RR$ are 
Lebesgue measurable functions, 
$g(t)\leq t$, $h(t)\leq t$, $t \geq 0$,  
${\displaystyle 
\limsup_{t\rightarrow\infty}(t-g(t))< \infty,
~\limsup_{t\rightarrow\infty}(t-h(t))< \infty}$.


Together with (\ref{1}) consider for each $t_0\geq 0$ an initial value problem
\begin{equation}
\label{4}
\ddot{x}(t) +a(t)\dot{x}(g(t))+b(t)x(h(t))=f(t), ~ t\geq t_0,
\end{equation}
\begin{equation}
\label{5}
x(t)=\varphi(t), \dot{x}(t)=\psi(t),~t<t_0, ~x(t_0)=x_0,~ \dot{x}(t_0)=x_0^{'}.
\end{equation}

We also assume that the following hypothesis holds

(a3) $f:[t_0,\infty)\rightarrow \RR $ is a Lebesgue measurable locally 
essentially bounded function, 
$\varphi:(-\infty,t_0)\rightarrow \RR, \psi:(-\infty,t_0)\rightarrow \RR$
are Borel measurable bounded functions.
\vspace{2mm}

\noindent
{\bf Definition.} A function $x:\RR \rightarrow \RR$ with locally 
absolutely continuous on $[t_0,\infty)$ derivative $\dot{x}$ is called 
{\em a solution} of  problem (\ref{4}), (\ref{5}) 
if it satisfies equation (\ref{4}) for almost every 
$t\in [t_0, \infty)$ and equalities (\ref{5}) for $t\leq t_0$.
\vspace{2mm}

\noindent
{\bf Definition}. For each $s\geq 0$, the solution $X(t,s)$ of the problem
\begin{equation}
\label{6}
\begin{array}{l}
{ \displaystyle \ddot{x}(t) +a(t)\dot{x}(g(t))+b(t)x(h(t))=0,~ t\geq s,}
\\ { \displaystyle 
x(t)=0, \dot{x}(t)=0,~t< s,~ x(s)=0,~ \dot{x}(s)=1 } 
\end{array}
\end{equation}
is called {\em the fundamental function} of  equation (\ref{1}).

We assume $X(t,s)=0,~ 0\leq t <s$.
\vspace{2mm}

Let  functions $x_1 $ and $x_2$ be the solutions of the following 
equation
$$
\ddot{x}(t) +a(t)\dot{x}(g(t))+b(t)x(h(t))=0,
~ t\geq t_0,~x(t)=0,~\dot{x}(t)=0,~ t<t_0,
$$
with initial values $x(t_0)=1,~\dot{x}(t_0)=0$ for $ x_1$ and 
$x(t_0)=0,~\dot{x}(t_0)=1$ for $x_2$, respectively.

By definition $x_2(t)=X(t,t_0)$.

\newtheorem{uess}{Lemma}
\newtheorem{guess}{Theorem}
\newtheorem{corollary}{Corollary}
\begin{uess}$\cite{AS}$
Let (a1)-(a3) hold. Then there exists one and only one
solution of  problem (\ref{4}), (\ref{5}) that can be presented in the form
\begin{equation}
\label{7}
x(t)=x_1(t)x_0+x_2(t)x_0^{'}+\int_{t_0}^t \!\! X(t,s)f(s)ds-
\int_{t_0}^t \!\!\!\! X(t,s)[a(s)\psi(g(s))+ b(s)\varphi(h(s))]ds,
\end{equation}
where $\varphi(h(s))=0$ if $h(s)>t_0$ and $\psi(g(s)))=0$ if  $g(s)>t_0$.
\end{uess}

\noindent
{\bf Definition}. Eq.~(\ref{1}) is {\em (uniformly) exponentially stable}, 
if there exist
$M>0$, $\mu>0$, such that  the solution of problem (\ref{1}),(\ref{5})
has the estimate
$$
|x(t)|\leq M~e^{-\mu (t-t_0)}\left[ |x(t_0)|+
\sup_{t<t_0} \left( |\varphi(t)|+|\psi(t)| \right) \right],~~t\geq t_ 0,
$$
where $M$ and $\mu$ do not depend on $t_ 0$.
\vspace{2mm}

\noindent
{\bf Definition}.
The fundamental function  $X(t,s)$ of (\ref{1}) {\em has
an exponential estimate}
if there exist
positive numbers 
$K>0, \lambda>0$, such that
\begin{equation}
\label{8}
|X(t,s)|\leq K~e^{-\lambda (t-s)},~~t\geq s\geq 0.
\end{equation}

For the linear equation (\ref{1})  with bounded delays ((a2) holds) 
the last  two definitions are equivalent. 

Under (a2)  the exponential 
stability does not depend on
values of equation parameters on any finite interval.
\vspace{2mm}

\noindent
{\bf Remark.}
All definitions and Lemma 1  can also be applied to equation 
(\ref{2}), for which we will assume that conditions (a2)-(a3) hold
 and all coefficients are Lebesgue 
measurable essentially bounded on $[0, \infty)$ functions.
\vspace{2mm}

Consider the equation
\begin{equation}
\label{9}
\ddot{x}(t) +a\dot{x}(t)+bx(t)=0, 
\end{equation}
where $a>0, b>0$ are positive numbers. This equation is exponentially stable.
Denote by $Y(t,s)$ the fundamental function of  (\ref{9}).

\begin{uess}$\cite{BBD}$
Let $a>0, b>0$.

1) If  $a^2>4b$ then ${\displaystyle \int_0^t|Y(t,s)|ds\leq\frac{1}{b}, 
~~ 
\int_0^t|Y^{'}_t(t,s)|ds\leq\frac{2a}{\sqrt{a^2-4b}(a-\sqrt{a^2-4b})} }$.
\vspace{2mm}

2) If $a^2<4b$ then ${\displaystyle \int_0^t|Y(t,s)|ds\leq 
\frac{4}{a\sqrt{4b-a^2}}, ~~
\int_0^t|Y^{'}_t(t,s)|ds\leq \frac{2(a+\sqrt{4b-a^2})}{a\sqrt{4b-a^2}} }$.
\vspace{2mm}

3)  If  $a^2=4b$ then ${\displaystyle \int_0^t|Y(t,s)|ds\leq \frac{1}{b},
~~\int_0^t|Y^{'}_t(t,s)|ds\leq \frac{2}{\sqrt{b}}}$.
\end{uess}

Let us introduce some functional spaces on a semi-axis. Denote
by ${\bf L}_{\infty}[t_0,\infty)$ the space of all essentially bounded on
$[t_0,\infty)$ scalar
functions and by ${\bf C}[t_0,\infty)$ the space of all continuous
bounded on  $[t_0,\infty)$ scalar functions
with the supremum norm. 

\begin{uess}$\cite{AS}$
Suppose there exists $t_0\geq 0$ such that for every 
$f\in {\bf L}_{\infty}[t_0,\infty)$ 
both the solution $x$ of the problem 
$$
\ddot{x}(t) +a(t)\dot{x}(g(t))+b(t)x(h(t))=f(t),~ t\geq t_0,
$$$$
x(t)=0, \dot{x}(t)=0,~t\leq t_0,
$$
and its derivative $\dot{x}$ belong to ${\bf C}[t_0,\infty)$. Then 
equation (\ref{1}) is exponentially
stable.
\end{uess}
{\bf Remark.}
A similar result is valid for equation (\ref{2}).

\begin{uess} $\cite{BB2}$
If $a(t)\geq 0$ is essentially bounded on $[0, \infty)$, the fundamental 
function $Z(t,s)$ of the equation
\begin{equation}
\label{10}
\dot{x}(t)+a(t)x(g(t))=0 
\end{equation} 
is positive: $Z(t,s)>0$, $t\geq s\geq t_0\geq 0$ and $t-g(t)\leq \delta$,
then  
$$
\int_{t_0+\delta}^t Z(t,s) a(s)ds\leq 1 \mbox{~~ for all ~~~} t \geq t_0+\delta.
$$ 
\end{uess}

\begin{uess} $\cite{BB2}$
\label{lemma5new}
If $a(t)\geq \alpha > 0$ is essentially bounded in $[0, \infty)$,
$\limsup_{t \to \infty} (t-g(t))<\infty$ and  the fundamental 
function $Z(t,s)$ of equation (\ref{10}) is positive, then 
(\ref{10}) is exponentially stable and $Z(t,s)$ has an exponential 
estimate.
\end{uess}

\begin{uess} $\cite{GL}$
Suppose $a(t)\geq 0$ and
$$
\int_{g(t)}^t a(s)ds\leq \frac{1}{e}, ~t\geq t_0\geq 0.
$$
Then the fundamental function of (\ref{10}) is positive:
$Z(t,s)>0, t\geq s\geq t_0$. 
\end{uess}

\section{Stability Conditions, I}

In this section we consider equation (\ref{2}) 
as a perturbation of an exponentially stable ordinary differential 
equation for which integral estimations of the fundamental function 
and its derivative are known.

We will start with the main result of \cite{BD}.

Denote by $Y(t,s)$ the fundamental function of the equation
\begin{equation}
\label{11}
\ddot{x}(t) +a(t)\dot{x}(t)+b(t)x(t)=0. 
\end{equation}
If equation (\ref{11}) is exponentially stable then there exist

\begin{equation}
\label{11a}
Y=\sup_{t>t_0}\int_{t_0}^t|Y(t,s)|ds<\infty, ~~
Y^{'}=\sup_{t>t_0}\int_{t_0}^t|Y^{'}_t(t,s)|ds<\infty.
\end{equation}
Denote
$$
a=\left\{\begin{array}{ll}
0,& g(t)\equiv t,\\ \sup_{t>t_0} a(t),&g(t)\not \equiv t,
\end{array}\right. ~~~~
b=\left\{\begin{array}{ll}
0,& h(t)\equiv t,\\ \sup_{t>t_0} b(t),&h(t)\not \equiv t,
\end{array}\right. 
$$$$
h=\max\left\{ \sup_{t>t_0} (t-g(t)),\sup_{t>t_0} (t-h(t)) \right\}.
$$
{\bf Theorem A} $\cite{BD}$
{\em Suppose $a(t)\geq 0$, $b(t)\geq 0$, equation (\ref{11}) is
exponentially 
stable and for some $t_0\geq 0$
$$
h<\frac{1}{a^2Y^{'}+abY+bY^{'}}.
$$
Then (\ref{1}) is exponentially stable.
}
\vspace{2mm} 

Now we will apply the method of \cite{BD} to equation (\ref{2}). 

\begin{guess}
Suppose $a(t)\geq 0, b(t)\geq 0$, equation (\ref{11}) is exponentially   
stable and for some $t_0\geq 0$
$$
\|a_1\|Y^{'}+ \|b_1\|Y<1,
$$
where $\|\cdot\|$ is the norm in the space ${\bf L}_{\infty}[t_0,\infty)$.
Then (\ref{2}) is exponentially stable.   
\end{guess}
{\bf Proof.}
Without loss of generality we can assume $t_0=0$. 
Consider the following problem 
\begin{equation}
\label{12}
\begin{array}{l}
\ddot{x}(t) +a(t)\dot{x}(t)+b(t)x(t)+a_1(t)\dot{x}(g(t))+b_1(t)x(h(t))=f(t),
~t \geq 0,\\
x(t)=\dot{x}(t)=0, ~t\leq 0.
\end{array}
\end{equation}

Let us demonstrate that for every $f\in {\bf L}_{\infty}$
the solution $x$ of (\ref{12}) and its derivative are bounded.
Due to the zero initial conditions, we can  assume that
$h(t)=g(t)=0,~t<0$.

The solution $x$ of (\ref{12}) is also a solution of the following 
problem
\begin{equation}
\label{13}
\ddot{x}(t) +a(t)\dot{x}(t)+b(t)x(t)=z(t), ~~t \geq 0, ~~ 
x(t)=\dot{x}(t)=0, t\leq 0,
\end{equation}
with some function $z(t)$. Then
\begin{equation}
\label{14}
x(t)=\int_0^t Y(t,s)z(s)ds,~ x'(t)=\int_0^t Y'_t(t,s)z(s)ds,
\end{equation}
where $Y(t,s)$ is the fundamental function of equation (\ref{11}).
Hence equation (\ref{12}) is equivalent to the equation
\begin{equation}
\label{15}
z(t)+a_1(t)\int_0^{g(t)}\!\!\!\!  Y'(g(t),s)z(s)ds+
b_1(t)\int_0^{h(t)}\!\!\!\!  Y(h(t),s)z(s)ds=f(t).
\end{equation}
Equation (\ref{15}) has the form $z+Hz=f$, which we consider in the space 
${\bf L}_{\infty}[0,\infty)$. 
We have
$$
\|H\|\leq \|a_1\|Y^{'}+ \|b_1\|Y<1.
$$
Then  for the solution of (\ref{15}) we have $z\in {\bf L}_{\infty}$.

Equalities (\ref{14}) imply $\|x\|\leq Y \|z\|,~ \|x'\|\leq Y' \|z\|$.

By Lemma 3 equation (\ref{2}) is exponentially stable.
\qed

Consider now equation (\ref{2}) with $a(t)\equiv a, b(t)\equiv b$.
By Lemma 2 we have the following statement.
\begin{corollary}
Suppose for some $ t_0\geq 0$ one of the following conditions holds:

1) $a^2>4b,$
$$
\frac{2a}{\sqrt{a^2-4b}(a-\sqrt{a^2-4b})}\|a_1\|+\frac{1}{b}\|b_1\|<1,
$$

2) $4b>a^2$,
$$
\frac{2(a+\sqrt{4b-a^2})}{a\sqrt{4b-a^2}}\|a_1\|+
\frac{4}{a\sqrt{4b-a^2}}\|b_1\|<1,
$$

3) $a^2=4b$,
$$
\frac{2}{\sqrt{b}}\|a_1\|+\frac{1}{b}\|b_1\|<1,
$$
where $\| \cdot \|$ is the norm in the space ${\bf 
L}_{\infty}[t_0,\infty)$.
Then equation (\ref{2}) is exponentially stable.
\end{corollary}

\begin{figure}[ht]
\centering
\setlength{\unitlength}{0.0005in}
\begingroup\makeatletter\ifx\SetFigFont\undefined%
\gdef\SetFigFont#1#2#3#4#5{%
  \reset@font\fontsize{#1}{#2pt}%
  \fontfamily{#3}\fontseries{#4}\fontshape{#5}%
  \selectfont}%
\fi\endgroup%
{\renewcommand{\dashlinestretch}{30}
\begin{picture}(4321,4539)(0,-10)
\path(2262,4512)(2337,4287)
\path(4062,2037)(4287,2112)
\path(2262,3912)(1962,2112)(2262,312)
	(2562,2112)(2262,3912)
\path(2412,3612)(2637,3837)
\path(2487,3762)(2412,3612)(2562,3687)
\path(2262,12)(2262,4512)(2187,4287)
\path(2862,1512)(2637,1737)(2712,1587)
\path(2787,1662)(2637,1737)
\path(12,2112)(4287,2112)(4062,2187)
\path(2637,2112)(2262,2787)(1887,2112)
	(2262,1437)(2637,2112)
\put(2112,1887){\makebox(0,0)[lb]{{\SetFigFont{10}{14.4}{\familydefault}{\mddefault}{\updefault}0}}}
\put(4137,2262){\makebox(0,0)[lb]{{\SetFigFont{10}{14.4}{\familydefault}{\mddefault}{\updefault}$a_1$}}}
\put(2412,4362){\makebox(0,0)[lb]{{\SetFigFont{10}{14.4}{\familydefault}{\mddefault}{\updefault}$b_1$}}}
\put(3237,1962){\makebox(0,0)[lb]{{\SetFigFont{10}{14.4}{\familydefault}{\mddefault}{\updefault}1}}}
\put(1287,1887){\makebox(0,0)[lb]{{\SetFigFont{10}{14.4}{\familydefault}{\mddefault}{\updefault}-1}}}
\put(2037,1212){\makebox(0,0)[lb]{{\SetFigFont{10}{14.4}{\familydefault}{\mddefault}{\updefault}-1}}}
\put(2187,3012){\makebox(0,0)[lb]{{\SetFigFont{10}{14.4}{\familydefault}{\mddefault}{\updefault}1}}}
\put(2037,312){\makebox(0,0)[lb]{{\SetFigFont{10}{14.4}{\familydefault}{\mddefault}{\updefault}-2}}}
\put(2712,3837){\makebox(0,0)[lb]{{\SetFigFont{10}{14.4}{\familydefault}{\mddefault}{\updefault}$a=3,b=2$}}}
\put(2937,1437){\makebox(0,0)[lb]{{\SetFigFont{10}{14.4}{\familydefault}{\mddefault}{\updefault}$a=3,b=2.5$}}}
\put(2112,3762){\makebox(0,0)[lb]{{\SetFigFont{10}{14.4}{\rmdefault}{\mddefault}{\updefault}2}}}
\end{picture}
}
\caption{The domain inside the rhombus gives the area of parameters
$a_1$, $b_1$ for the exponential stability of (\protect{\ref{2}}): 
for $a=3$, $b=2$ and $a=3$, $b=2.5$, respectively.}
\label{figure1}
\end{figure}
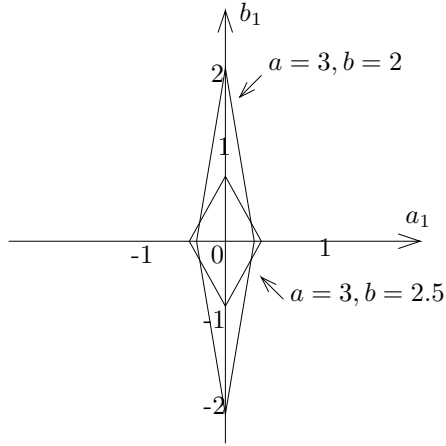

{\bf Example 1.} Let us illustrate the exponential stability domain 
in Corollary 1 for two cases. If $a=3$, $b=2$ then $a^2>4b$ and the 
condition of 1) becomes ${\displaystyle 3|a_1|+0.5 |b_1| <1}$
for constant $a_1$, $b_1$ which corresponds to the domain inside the 
vertically stretched rhombus in Fig.~\ref{figure1}. If $a=3$, $b=2.5$ 
then $a^2<4b$ and the condition of 2) becomes ${\displaystyle 2|a_1|+|b_1| 
< 0.75}$ which is inside a smaller rhombus in Fig.~\ref{figure1}.

\section{Stability Conditions II}

In this section we consider equations (\ref{1}) and (\ref{2}) 
as perturbations of an exponentially stable ordinary differential 
equation for which only an integral estimation of 
its fundamental function is known.

\begin{guess}
Suppose  (\ref{11}) is exponentially stable, (\ref{10}) has 
a positive fundamental function  $Z(t,s)>0$, $t\geq s\geq t_0\geq 0$,  
~$a(t)\geq \alpha>0$, ~$t-g(t)\leq \delta$, $t-h(t)\leq \tau$.
 
If for some $t_1\geq t_0+\delta$ 
\begin{equation}
\label{18}
Y\left[\delta\|a\|\left(\|a\|\left\|\frac{b}{a}\right\|+\|b\|\right)
+\tau\|b\|\left\|\frac{b}{a}\right\|\right]<1,
\end{equation}
where $Y$ is denoted in (\ref{11a}), $\| \cdot \|$ is the norm in the 
space ${\bf L}_{\infty}[t_1,\infty)$, 
then  equation (\ref{1})
is exponentially stable.
\end{guess}
{\bf Proof.}
Consider the following problem 
\begin{equation}
\label{19}
\begin{array}{l}
\ddot{x}(t) +a(t)\dot{x}(g(t))+b(t)x(h(t))=f(t),\\
x(t)=\dot{x}(t)=0, ~t\leq t_1.
\end{array}
\end{equation}
As in the proof of Theorem 1 we can 
assume that $h(t)=g(t)=0$, $t<t_1$.
Let us demonstrate that for every $f\in {\bf L}_{\infty}[t_1,\infty)$
the solution of problem (\ref{19}) and its derivative are bounded.

The equation in problem (\ref{19}) can be rewritten as
\begin{equation}
\label{20}   
\ddot{x}(t) +a(t)\dot{x}(t)+b(t)x(t)-a(t)\int_{g(t)}^t \ddot{x}(s)ds
-b(t)\int_{h(t)}^t \dot{x}(s)ds=f(t).
\end{equation}
Equation (\ref{20}) is equivalent to the following one
\begin{equation}
\label{21}
x(t)-\int_{t_1}^t Y(t,s)a(s)\int_{g(s)}^s \ddot{x}(\tau)d\tau\, ds 
-\int_{t_1}^t Y(t,s)b(s)\int_{h(s)}^s \dot{x}(\tau)d\tau\, ds =f_1(t),
\end{equation}
where $f_1(t)=\int_{t_1}^t Y(t,s)f(s)ds$ and $Y(t,s)$ is the fundamental 
function of equation (\ref{11}). Hence $f_1\in {\bf 
L}_{\infty}[t_1,\infty)$.

The equation in problem (\ref{19}) can be rewritten in a different form
\begin{equation}
\label{22}
\dot{x}(t)+\int_{t_1}^t Z(t,s)b(s)x(h(s))ds=r(t),
\end{equation}
where $r(t)=\int_{t_1}^t Z(t,s)f(s)ds$, $Z(t,s)$ is the fundamental 
function of (\ref{10}). Since $a(t)\geq \alpha>0$ 
and $Z(t,s)>0$, then by Lemma \ref{lemma5new} fundamental function 
$Z(t,s)$ has an exponential estimation. Hence
$r\in {\bf L}_{\infty}[t_1,\infty)$.

From (\ref{19}) and (\ref{21}) we have
\begin{eqnarray}
x(t) & + &\int_{t_1}^t Y(t,s)a(s)\int_{g(s)}^s
\left[ a(\tau)\dot{x}(g(\tau))+b(\tau)x(h(\tau))\right] d\tau \, ds
\nonumber \\
\label{23}
& - & \int_{t_1}^t Y(t,s)b(s)\int_{h(s)}^s\dot{x}(\tau)d\tau ds=f_2(t),
\end{eqnarray}
where 
$$
f_2(t)=f_1(t)+\int_{t_1}^t Y(t,s)a(s)\int_{g(s)}^s f(\tau)d\tau\, ds.
$$
Since $\|f_2\|\leq \|f_1\|+Y\delta \|a\| \|f\|$, then
$f_2\in {\bf L}_{\infty}[t_1,\infty)$.

Substituting $\dot{x}$  from (\ref{22}) into (\ref{23}), we obtain
$$
x(t)  - \int_{t_1}^t Y(t,s)a(s)\int_{g(s)}^s
\left[a(\tau)\int_{t_1}^{g(\tau)}
Z(g(\tau),\xi)b(\xi)x(h(\xi))d\xi -b(\tau)x(h(\tau))\right]d\tau\, ds
$$
\begin{equation}
\label{24}
+\int_{t_1}^t Y(t,s)b(s)\int_{h(s)}^s \left[ \int_{t_1}^{\tau} 
Z(\tau,\xi)b(\xi)x(h(\xi))d\xi \right]
~d\tau ~ ds=f_3(t),
\end{equation}  
where 
$$
f_3(t)=f_2(t)-\int_{t_1}^t Y(t,s)\left[a(s)\int_{g(s)}^s
r(g(\tau))d\tau-b(s)\int_{h(s)}^sr(\tau)d\tau\right] ds.
$$
Since $\| f_3\| \leq \|f_2 \|+Y\| r \| (\|a\|\delta+ \|b\| \tau)$
then $f_3\in {\bf L}_{\infty}[t_1,\infty)$.

Equation (\ref{24}) has the form $x-Tx=f_3$. 
Lemma 4 yields that
$$
\int_{t_1}^{g(t)} Z(g(t),s)b(s)ds\leq \sup_{t>t_1}
\int_{t_1}^t Z(t,s)b(s)ds\leq \sup_{t>t_1} \int_{t_1}^t Z(t,s) a(s)\frac{b(s)}{a(s)}ds
\leq \left\|\frac{b}{a}\right\|.
$$


Condition  (\ref{18}) implies $\|T\|<1$.
Hence $x\in {\bf C}[t_1,\infty)$. From (\ref{22}) we have $\dot{x}\in 
{\bf L}_{\infty}[t_1,\infty)$.
By Lemma 3 equation (\ref{1}) is exponentially stable.
\qed

\noindent
{\bf Remark.}
Lemma 6 gives explicit conditions for positivity of 
the fundamental function of equation 
(\ref{10}). 
\vspace{2mm}

Consider the equation with constant coefficients
\begin{equation}
\label{25}
\ddot{x}(t) +a\dot{x}(g(t))+bx(h(t))=0,
\end{equation}
where $a>0, b>0, t-g(t)\leq \delta, t-h(t)\leq \tau$.
\begin{corollary}
Suppose $a\delta\leq \frac{1}{e}$ and one of the following conditions 
holds:

1) ${\displaystyle a^2\geq 4b, ~ 2\delta a+\frac{\tau b}{a}<1,}$
\vspace{2mm}

2) ${\displaystyle a^2<4b,~2\delta ab+\frac{\tau 
b^2}{a}<\frac{a\sqrt{4b-a^2}}{4}}$.
\vspace{2mm}

Then equation (\ref{25}) is exponentially stable.
\end{corollary}

\noindent
{\bf Example 2.} Let us illustrate the exponential stability domain 
in Corollary 2 for both cases. If $a=3$, $b=2$ then $a^2>4b$ and the 
condition of 1) becomes ${\displaystyle 6\delta+\frac{2}{3} \tau 
<1}$. If $a=3$, $b=2.5$ then $a^2<4b$ and the
condition of 2) becomes ${\displaystyle 15 \delta+ 
\frac{25}{12}\tau<\frac{3}{4} }$, or ${\displaystyle 
20\delta+\frac{25}{9}\tau<1}$. Here the 
inequality $a\delta=3\delta 
<1/e$ should also be satisfied (the area under the horizontal line).

\begin{figure}[ht]
\centering
\setlength{\unitlength}{0.0007in}
\begingroup\makeatletter\ifx\SetFigFont\undefined%
\gdef\SetFigFont#1#2#3#4#5{%
  \reset@font\fontsize{#1}{#2pt}%
  \fontfamily{#3}\fontseries{#4}\fontshape{#5}%
  \selectfont}%
\fi\endgroup%
{\renewcommand{\dashlinestretch}{30}
\begin{picture}(6483,2139)(0,-10)
\path(87,612)(6462,612)(6237,687)
\path(6237,537)(6462,612)
\path(12,987)(5862,987)
\path(687,1212)(6087,612)
\path(2487,1587)(2187,1062)(2262,1287)
\path(2187,1062)(2337,1212)
\path(687,12)(687,2112)(612,1812)
\path(687,2112)(762,1812)
\path(687,762)(2112,612)
\path(1137,1662)(912,762)(912,912)
\path(912,762)(987,912)
\put(6237,762){\makebox(0,0)[lb]{\smash{{{\SetFigFont{10}{14.4}{\familydefault}{\mddefault}{\updefault}$\tau$}}}}}
\put(1962,387){\makebox(0,0)[lb]{\smash{{{\SetFigFont{10}{14.4}{\familydefault}{\mddefault}{\updefault}$\frac{9}{25}$}}}}}
\put(5937,387){\makebox(0,0)[lb]{\smash{{{\SetFigFont{10}{14.4}{\familydefault}{\mddefault}{\updefault}$\frac{3}{2}$}}}}}
\put(537,387){\makebox(0,0)[lb]{\smash{{{\SetFigFont{10}{14.4}{\familydefault}{\mddefault}{\updefault}0}}}}}
\put(462,1137){\makebox(0,0)[lb]{\smash{{{\SetFigFont{10}{14.4}{\familydefault}{\mddefault}{\updefault}$\frac{1}{6}$}}}}}
\put(2562,1512){\makebox(0,0)[lb]{\smash{{{\SetFigFont{10}{14.4}{\familydefault}{\mddefault}{\updefault}$a=3,b=2$}}}}}
\put(462,1812){\makebox(0,0)[lb]{\smash{{{\SetFigFont{10}{14.4}{\familydefault}{\mddefault}{\updefault}$\delta$}}}}}
\put(1062,1662){\makebox(0,0)[lb]{\smash{{{\SetFigFont{10}{14.4}{\familydefault}{\mddefault}{\updefault}$a=3,b=2.5$}}}}}
\put(5037,1062){\makebox(0,0)[lb]{\smash{{{\SetFigFont{10}{14.4}{\familydefault}{\mddefault}{\updefault}$a 
\delta=\frac{1}{e}$}}}}}
\put(462,687){\makebox(0,0)[lb]{\smash{{{\SetFigFont{10}{14.4}{\familydefault}{\mddefault}{\updefault}$\frac{1}{20}$}}}}}
\end{picture}
}
\caption{The domain inside the triangle under the horizontal line
$a\sigma=3\sigma= 1/e$ gives the domain of parameters
$\tau$, $\sigma$ for the exponential stability of 
(\protect{\ref{25}}) : 
for $a=3$, $b=2$ and $a=3$, $b=2.5$, respectively. The domain of 
parameters in the former case involves the domain in the latter case.}
\label{figure2}
\end{figure}
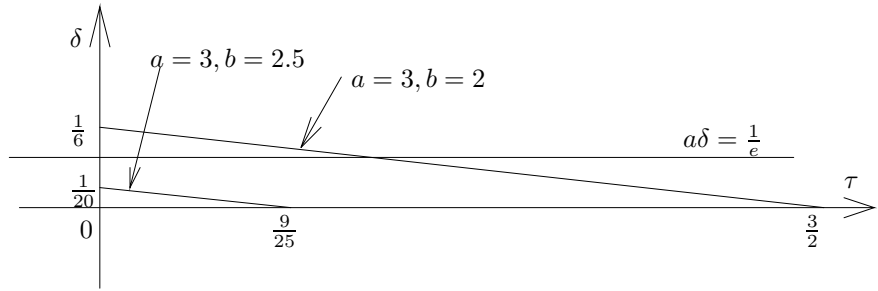

\begin{corollary}
Suppose $g(t)\equiv t, a^2\geq 4b, \tau b<a$.
Then equation (\ref{25}) is exponentially stable. 
\end{corollary}

\noindent
{\bf Remark.}
Corollary 3 gives the same condition $\tau b<a$ as was obtained 
by Burton in \cite{Bur} for autonomous equation  (\ref{3});
however, equation (\ref{25}) is not autonomous.
\vspace{2mm}

In Theorem 2 it was assumed that the first order equation (\ref{10})
has a positive fundamental function.
In the next theorem we will omit this restriction.
\begin{guess}
Suppose  (\ref{11}) is exponentially stable, 
$a(t)\geq \alpha>0$, $t-g(t)\leq \delta$, $t-h(t)\leq \tau$.

If for some $t_0\geq 0$ the inequality $\delta\|a\|<1$ holds and
\begin{equation}
\label{26}
Y\left[\frac{(\delta\|a\|^2+\tau\|b\|)(\left\|\frac{b}{a}\right\|+\delta 
\|b\|)}{1-\delta \|a\|}+\delta \|b \|\right]<1,
\end{equation}
where $\|  \cdot \|$ is the norm in the space 
${\bf L}_{\infty}[t_0,\infty)$,
then  equation (\ref{1}) is exponentially stable.
\end{guess}
{\bf Proof.}
Without loss of generality we can assume $t_0=0$ and $h(t)=g(t)=0, t<0$.
As in the proof of Theorem 2 
we will demonstrate that for every $f\in {\bf L}_{\infty}[0,\infty)$
the solution of (\ref{19}) and its derivative are bounded.

First we will obtain an apriori estimation for the derivative 
of the solution of (\ref{19}).
Equation (\ref{19}) can be rewritten in the form
\begin{equation}
\label{27}
\ddot{x}(t) +a(t)\dot{x}(t)-a(t)\int_{g(t)}^t \ddot{x}(s)ds+b(t)x(h(t))=f(t).
\end{equation}
Substituting $\ddot{x}$ from (\ref{19}) into (\ref{27}) we have
\begin{equation}
\label{28}
\ddot{x}(t) +a(t)\dot{x}(t)+a(t)\int_{g(t)}^t 
[a(s)\dot{x}(g(s))+b(s)x(h(s))]ds
+b(t)x(h(t))=r(t),
\end{equation}
where $r(t)=f(t)+a(t)\int_{g(t)}^t f(s)ds$. 
Evidently $r\in {\bf L}_{\infty}[0,\infty)$.

Equation (\ref{28}) is equivalent to the following one
\begin{equation}
\label{29}
\dot{x}(t)+\int_0^t \!\! e^{-\int_s^t a(\xi)d\xi} 
\left(a(s)\int_{g(s)}^s \!\!        
\left[ 
a(\tau)\dot{x}(g(\tau))+b(\tau)x(h(\tau))\right]d\tau+b(s)x(h(s))\right) 
ds=r_1(t),
\end{equation}
where $r_1(t)=\int_0^t e^{-\int_s^t a(\xi)d\xi}r(s)ds$. 
Since $a(t)\geq\alpha>0$, then $r_1\in {\bf L}_{\infty}[0,\infty)$.

Denote by $\| \cdot \|_T$ the norm in the space ${\bf L}_{\infty}[0,T]$.
From (\ref{29}) we have 
$$
\|\dot{x} \|_T\leq
\delta( \|a \| \| \dot{x}\|_T+ 
\|b \| \|x \|_T)+\left\|\frac{b}{a}\right\| \|x \|_T+\|r_1 \|.
$$
Hence we obtain the following apriori estimation for $\| \dot{x} \|_T$
\begin{equation}
\label{30}
\|\dot{x} \|_T\leq \frac{\delta \|b\|+\left\|\frac{b}{a}\right\|}{ 
1-\delta \| a \|}\|x\|_T
+\frac{\| r_1 \| }{ 1-\delta \|a \|}.
\end{equation}
Equation (\ref{19}) is equivalent to (\ref{21}) which can be rewritten as
\begin{eqnarray}
x(t)& + & \int_0^t  Y(t,s)a(s)\int_{g(s)}^s 
[a(\tau)\dot{x}(g(\tau))+b(\tau)x(h(\tau)) ]d\tau
\, ds \nonumber
\\
\label{31}
& - & \int_0^t \! Y(t,s)b(s)\int_{h(s)}^s \dot{x}(\tau)
d\tau \, ds  = f_2(t),
\end{eqnarray}
where $f_2(t)=f_1(t)+\int_0^t Y(t,s)a(s)\int_{g(s)}^s f(\tau)d\tau ds$,
$Y(t,s)$ is the fundamental function of (\ref{11}). 
Then $f_2\in {\bf L}_{\infty}[0,\infty)$ and
\begin{equation}
\label{32}
\|x\|_T\leq Y \left[ \delta \|a\|(\| a \| \| \dot{x}\|_T+\|b \| \|x\|_T)+
\tau \| b \| \| \dot{x}\|_T \right] + \| f_2 \|
\end{equation}
$$
=Y[(\left[ \delta \| a \|^2+\tau 
\| b \|) \| \dot{x} \|_T+\delta \|b\| \|x\|_T \right] + \|f_2\|.
$$
Inequalities  (\ref{30}) and (\ref{32}) imply
\begin{equation}
\label{33}
\| x \|_T\leq Y\left[\frac{(\delta \|a\|^2+\tau \|b \|)
(\left\| \frac{b}{a}\right\|+\delta \| b \|)}   
{1-\delta \|a \|}+\delta \| b \|\right] \| x \|_T +C, 
\end{equation}
where $C$ is a  positive number.
 
Inequality (\ref{33}) has the form $\|x \|_T\leq A\|x\|_T+C$, where $A$ 
and 
$C$ do not depend on $T$ and $0<A<1$. Hence $x\in 
{\bf L}_{\infty}[0,\infty)$.
Inequality (\ref{30}) implies  $\dot{x}\in {\bf L}_{\infty}[0,\infty)$.
By Lemma 3 equation (\ref{1}) is exponentially stable.
\qed 

Now let us apply the method of apriori estimation to equation (\ref{2}).
\begin{guess}
Suppose  (\ref{11}) is exponentially stable,
$a(t)\geq \alpha>0 $, $t-g(t)\leq \delta$, $t-h(t)\leq \tau$.

If for some $t_0\geq 0$ the inequality $\left\|\frac{a_1}{a}\right\|<1$ 
holds and
\begin{equation}
\label{34}
Y\left[\frac{\|a_1\|\left(\left\|\frac{b}{a}\right\|
+\left\|\frac{b_1}{a}\right\|\right)}
{1-\left\|\frac{a_1}{a}\right\|}+\|b_1 \|\right]<1,
\end{equation}
where $\| \cdot \|$ is the norm in the space ${\bf L}_{\infty}[t_0,\infty)$,
$Y$ is defined in (\ref{11a}), then equation (\ref{2})
is exponentially stable.
\end{guess}
{\bf Proof.}
Without loss of generality we can assume $t_0=0$.
Let us demonstrate that for every $f\in {\bf L}_{\infty}[0,\infty)$
both the solution $x$ of (\ref{12}) and its derivative are bounded.

Equation (\ref{12}) is equivalent to the following one
\begin{equation}
\label{35} 
\dot{x}(t)=-\int_0^t e^{-\int_s^t a(\tau) d\tau} 
[a_1(s)\dot{x}(g(s))+b_1(s)x(h(s))+b(s)x(s)]ds+r(t),
\end{equation}
where $r(t)=\int_0^t e^{-\int_s^t a(\tau) d\tau}f(s)ds$.
Evidently $r\in {\bf L}_{\infty}[0,\infty)$. Hence 
$$
\|\dot{x}\|_T\leq \left\|\frac{a_1}{a}\right\| \|\dot{x}\|_T+
\left(\left\|\frac{b}{a}\right\|+\left\|\frac{b_1}{a}\right\|\right)
\|x \|_T+ \|r \|,
$$
where $\| \cdot \|_T$ is the norm in ${\bf L}_{\infty}[0,T]$.
Then
\begin{equation}
\label{36}
\| \dot{x} \|_T \leq 
\frac{\left\| \frac{b}{a} \right\| + \left\|\frac{b_1}{a}\right\|}
{1-\left\|\frac{a_1}{a}\right\|}\| x \|_T+\frac{\|r\|}{1-\left\|\frac{a_1}{a}
\right\|}.
\end{equation}
Equation (\ref{12}) is also equivalent to 
\begin{equation}
\label{37}
x(t)=-\int_0^t Y(t,s)\left[a_1(s)\dot{x}(g(s))+b_1(s)x(h(s))\right] 
ds+f_1(t),
\end{equation}
where $f_1(t)=\int_0^t Y(t,s) f(s)ds$,
$Y(t,s)$ is the fundamental function of (\ref{11}). The first inequality 
in (\ref{11a})
implies $f_1\in {\bf L}_{\infty}[0,\infty)$.
Thus
$$
\|x\|_T\leq Y(\| a_1 \| \| \dot{x} \|_T+\|b_1\| \|x\|_T)+\|f_1\|,
$$
which together with inequality (\ref{36}) yields
$$
\| x \|_T\leq 
Y\left[\frac{\| 
a_1\|\left(\left\|\frac{b}{a}\right\|+\left\|\frac{b_1}{a}\right\|\right)}
{1-\left\|\frac{a_1}{a}\right\|}+ \| b_1 \| \right] \| x \|_T +C.
$$
Here $C$ is a positive number which does not depend on $x$ and $T$.

As in the proof of Theorem 3, we obtain that $x$ and $\dot{x}$ 
are bounded functions.
Then by Lemma 3 equation (\ref{2}) is exponentially stable.
\qed

\noindent
{\bf Remark.} If $a_1(t)\equiv 0$ then the statements of 
Theorems 1 and 4 coincide.

\begin{corollary}
Suppose $a(t)\equiv a>0$, $b(t)\equiv b>0$, $a^2\geq 4b, t-g(t)\leq 
\delta$, $t-h(t)\leq \tau$ 
and for some $t_0\geq 0$ we have $a> \| a_1 \|,~\|a_1\|(b+\|b_1\|)<(a-\|a_1\|)(b-\|b_1\|)$,
where $\|\cdot \|$ is the norm in the space ${\bf 
L}_{\infty}[t_0,\infty)$.
Then equation (\ref{2}) is exponentially stable. 
\end{corollary}

\begin{figure}[ht]
\centering
\setlength{\unitlength}{0.0004in}
\begingroup\makeatletter\ifx\SetFigFont\undefined%
\gdef\SetFigFont#1#2#3#4#5{%
  \reset@font\fontsize{#1}{#2pt}%
  \fontfamily{#3}\fontseries{#4}\fontshape{#5}%
  \selectfont}%
\fi\endgroup%
{\renewcommand{\dashlinestretch}{30}
\begin{picture}(4321,4539)(0,-10)
\path(2262,4512)(2337,4287)
\path(12,2112)(4287,2112)(4062,2187)
\path(4062,2037)(4287,2112)
\path(2262,12)(2262,4512)(2187,4287)
\path(3612,2112)(2262,3912)(912,2112)
	(2262,312)(3612,2112)
\put(2112,1887){\makebox(0,0)[lb]{\smash{{{\SetFigFont{10}{14.4}{\familydefault}{\mddefault}{\updefault}0}}}}}
\put(4137,2262){\makebox(0,0)[lb]{\smash{{{\SetFigFont{10}{14.4}{\familydefault}{\mddefault}{\updefault}$a_1$}}}}}
\put(2412,4362){\makebox(0,0)[lb]{\smash{{{\SetFigFont{10}{14.4}{\familydefault}{\mddefault}{\updefault}$b_1$}}}}}
\put(2112,3987){\makebox(0,0)[lb]{\smash{{{\SetFigFont{10}{14.4}{\familydefault}{\mddefault}{\updefault}2}}}}}
\put(3237,1962){\makebox(0,0)[lb]{\smash{{{\SetFigFont{10}{14.4}{\familydefault}{\mddefault}{\updefault}1}}}}}
\put(1287,1887){\makebox(0,0)[lb]{\smash{{{\SetFigFont{10}{14.4}{\familydefault}{\mddefault}{\updefault}-1}}}}}
\put(2037,1212){\makebox(0,0)[lb]{\smash{{{\SetFigFont{10}{14.4}{\familydefault}{\mddefault}{\updefault}-1}}}}}
\put(2187,3012){\makebox(0,0)[lb]{\smash{{{\SetFigFont{10}{14.4}{\familydefault}{\mddefault}{\updefault}1}}}}}
\put(2037,312){\makebox(0,0)[lb]{\smash{{{\SetFigFont{10}{14.4}{\familydefault}{\mddefault}{\updefault}-2}}}}}
\end{picture}
}
\caption{The domain inside a rhombus gives the area of parameters
$a_1$, $b_1$ for the exponential stability of (\protect{\ref{2}})  
for $a=3$, $b=2$.}
\label{figure3}
\end{figure}
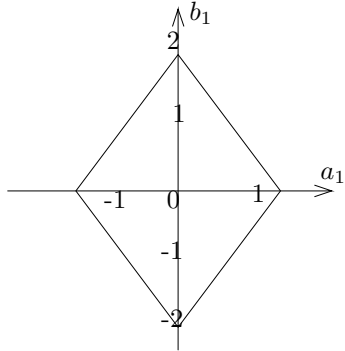

\noindent
{\bf Example 3.} Let us illustrate the exponential stability domain 
in Corollary 4. If $a=3$, $b=2$ then $a^2>4b$ and the 
condition becomes ${\displaystyle 4|a_1|+3 |b_1| <6}$
for constant $a_1$, $b_1$ (the domain inside the rhombus in 
Fig.~\ref{figure3}) and bounded delays.

\begin{corollary}
Suppose $a_1(t)\equiv 0$, $a(t)\equiv a>0$, $b(t)\equiv b>0$, $a^2\geq 
4b$, $t-h(t)\leq \tau$ 
and for some $t_0\geq 0$  we have  $\|b_1\|<b$,
where $\|\cdot \|$ is the norm in the space ${\bf L}_{\infty}[t_0,\infty)$.
Then equation (\ref{2}) is exponentially stable. 
\end{corollary}

\section{Explicit Stability Conditions }

In all previous results we have assumed that the ordinary differential 
equation (\ref{11}) is exponentially stable and an estimation 
on the integral of its fundamental function is known. 
All these conditions can hold, as the following lemma demonstrates.

\begin{uess}$\cite{BBD}$
Suppose for some $t_0\geq 0$
\begin{equation}
\label{38a}
a=\inf_{t\geq t_0}a(t)>0,~ b=\inf_{t\geq t_0}b(t)>0,~ B=\sup_{t\geq
t_0}b(t),~
a^2\geq 4B.
\end{equation}
Then the fundamental function of (\ref{11}) is positive: $Y(t,s)>0, 
t>s\geq t_0$.

Moreover, (\ref{11}) is exponentially stable and for its fundamental 
function we have
\begin{equation}
\label{38}
\int_{t_0}^t Y(t,s) b(s)ds \leq 1.
\end{equation}
\end{uess}
As an application of Lemma 7,  we will obtain  new explicit 
stability conditions for equations (\ref{1})  and (\ref{2}).

\begin{guess}
Suppose there 
exists  $t_0 \geq 0$ such that condition (\ref{38a}) holds, 
\begin{equation}
\label{39a}
\int_{g(t)}^t a(s)ds\leq \frac{1}{e}, ~t\geq t_0,
\end{equation}
and $t-g(t)\leq \delta$, $t-h(t)\leq \tau$ for $t\geq t_0$. 

If for some $t_1\geq t_0+\delta$
\begin{equation}
\label{39}
\delta\left\|\frac{a}{b}\right\|\left(\|a \|\left\|\frac{b}{a}\right\|+
\|b\| \right)+\tau\left\|\frac{b}{a}\right\|<1,
\end{equation}
where $\| \cdot \|$ is the norm in the space ${\bf 
L}_{\infty}[t_1,\infty)$, 
then equation (\ref{1})
is exponentially stable.
\end{guess}
{\bf Proof.}
 Let us demonstrate that for every $f\in {\bf L}_{\infty}[t_1,\infty)$
the solution $x$ of problem (\ref{19}) and its derivative are bounded.
The equation in problem (\ref{19})  is equivalent to 
(\ref{24}), 
which has the form $x-Hx=f_3$. Condition (\ref{39a}) 
and  Lemma 6 imply that the fundamental function of equation  (\ref{10})
satisfies $Z(t,s)>0, t\geq t_0\geq 0$. Lemmas 4 and 6 yield that 
\begin{eqnarray*}
|(Hx)(t)| & \leq & \int_{t_1}^t Y(t,s)b(s)\frac{a(s)}{b(s)}\int_{g(s)}^s
\left[a(\zeta)\int_{t_1}^{g(\zeta)}
Z(g(\zeta),\xi)a(\xi)\frac{b(\xi)}{a(\xi)}d\xi +b(\zeta)\right]d\zeta\, 
ds~\| x \|
\\
& + & \int_{t_1}^t Y(t,s)b(s)\left[ \int_{h(s)}^s d\zeta \int_{t_1}^\zeta 
Z(\zeta,\xi)a(\xi)\frac{b(\xi)}{a(\xi)}d\xi \right] ds\| x \|
\\
& \leq & 
\left[\delta\left\|\frac{a}{b}\right\|\left(\|a\|\left\|\frac{b}{a}\right\|
+\|b\|\right)+\tau\left\|\frac{b}{a}\right\|\right] \|x\|.
\end{eqnarray*}
Inequality (\ref{39}) implies $\|H\|<1$, hence the 
solution $x$ of (\ref{24})
and therefore of (\ref{1}) is a bounded  function.
Equality (\ref{22}) yields that $\dot{x}$ is a bounded function.
By Lemma 3 equation (\ref{1}) is exponentially stable.
\qed

\begin{corollary}
Suppose there 
exists  $t_0 \geq 0$ such that condition  (\ref{38a})
 holds , $g(t)\equiv t, t-h(t)\leq \tau, t\geq t_0$. If
$\tau\left\|\frac{b}{a}\right\|<1$, where $\| \cdot \|$ 
is the norm in the space ${\bf 
L}_{\infty}[t_0,\infty)$, then equation (\ref{1}) is exponentially stable.
\end{corollary}

Let us proceed to stability conditions for equation (\ref{2}).
\begin{guess}
Suppose there 
exists  $t_0 \geq 0$ such that condition  (\ref{38a})
 holds,   $t-g(t)\leq 
\delta, t-h(t)\leq \tau, t\geq t_0$. 

If $\left\|\frac{a_1}{a}\right\|<1$ and
$$
\left\|\frac{a_1}{b}\right\|\frac{\left\|\frac{b}{a}\right\|+
\left\|\frac{b_1}{a}\right\|}
{1-\left\|\frac{a_1}{a}\right\|}+\left\|\frac{b_1}{b}\right\|<1,
$$
where $\| \cdot \|$ is the norm in the space ${\bf 
L}_{\infty}[t_0,\infty)$,
then equation (\ref{2})
is exponentially stable.
\end{guess}
{\bf Proof.}
We follow the proof of Theorem 4. From  (\ref{36}) and (\ref{37})
we have
$$
\|x\|_T\leq \int_0^t Y(t,s)b(s)ds 
\left(\left\|\frac{a_1}{b}\right\| \| \dot{x} \|_T
+\left\|\frac{b_1}{b}\right\| \|x|\|_T\right)+\|f_1\|
$$$$
\leq
\left(\left\|\frac{a_1}{b}\right\|\frac{\left\|\frac{b}{a}\right\|+\left\|\frac{b_1}{a}\right\|}
{1-\left\|\frac{a_1}{a}\right\|}+\left\|\frac{b_1}{b}\right\|\right)\|x\|_T 
+C.
$$
The end of the proof is similar to the proof of Theorem 4.
\qed

\begin{corollary}
Suppose there 
exists  $t_0 \geq 0$ such that condition (\ref{38a})
 holds,   $a_1(t)\equiv 0, t-h(t)\leq \tau, t\geq t_0$. 

If $\left\|\frac{b_1}{b}\right\|<1$, where $\| \cdot \|$ is the norm in 
the space ${\bf L}_{\infty}[t_0,\infty)$, then  equation (\ref{2})
is exponentially stable .
\end{corollary}


\begin{thebibliography}{99}

\bibitem{AS}
N.\,V. Azbelev and P.\,M. Simonov, Stability of Differential
Equations with Aftereffect. {\em Stability and Control:
Theory, Methods and Applications}, {\bf 20}. Taylor $\&$ Francis, London,
2003. 

\bibitem{BD}
D. Bainov and A. Domoshnitsky,
Stability of a second-order differential equation with retarded argument,
{\em Dynamics and Stability of Systems} {\bf 9} (1994), no. 2, 145--151. 

\bibitem{BB2}
L. Berezansky and E. Braverman,
On exponential stability of linear differential equations with several
delays, {\em J. Math. Anal. Appl.} {\bf 324} (2006), no. 2, 1336--1355.

\bibitem{BBD}
L. Berezansky, E. Braverman and A. Domoshnitsky,
Positive solutions and stability of a linear ordinary differential 
equation of the second order with damping term,
arXiv:0807.2227v1 [math.DS] (July 14, 2008); to appear in
{\em Funct. Differ. Equ.}, 2009.
\bibitem{Bur}
T.\,A. Burton,
Stability and Periodic Solutions of Ordinary and Functional 
Differential Equations, Academic Press, {\em Mathematics in Science and 
Engineering} {\bf 178}, 1985.

\bibitem{Bur2}
T.\,A. Burton,
Stability by Fixed Point Theory for Functional Differential Equations.
Dover Publications, Mineola, New York, 2006.

\bibitem{Bur3}
T.\,A. Burton,
Fixed points, stability, and exact linearization,
{\em Nonlinear Anal.} {\bf 61} (2005), 857--870.

\bibitem{BurFur}  
T.\,A. Burton and T. Furumochi,
Asymptotic behavior of solutions of functional 
differential equations by fixed point theorems,
{\em Dynam. Systems Appl.} {\bf 11} (2002), no. 4, 499--519.

\bibitem{BH}
T.\,A. Burton and L. Hatvani,
Asymptotic stability of second order ordinary, 
functional, and partial differential equations,
{\em J. Math. Anal. Appl.} {\bf  176} (1993), 261--281.

\bibitem{CS}
B. Cahlon and D. Schmidt,
Stability criteria for certain second-order delay differential equations
with mixed coefficients,
{\em J. Comput. Appl. Math.} {\bf 170} (1994), no. 1, 79--102.

\bibitem{D}
A. Domoshnitsky,  Componentwise applicability of Chaplygin's theorem 
to a system of linear differential equations with delay, 
{\em Differential Equations} {\bf 26} (1991), no. 10, 1254--1259.

\bibitem{EKZ}
L.\,N. Erbe, Q. Kong and B.\,G. Zhang,
Oscillation Theory for Functional Differential Equations,
Marcel Dekker, New York, Basel, 1995.

\bibitem{GL}
I. Gy\"{o}ri and G. Ladas, 
Oscillation Theory of Delay Differential Equations,
Clarendon Press, Oxford, 1991.

\bibitem{KN}
V.\,B. Kolmanovsky and V.\,R. Nosov,
Stability of Functional Differential Equations.
Academic Press, 1986.

\bibitem{LLZ}
G.\,S. Ladde, V. Lakshmikantham  and B.\,G. Zhang ,
Oscillation Theory of Differential Equations 
with Deviating Argument,
Marcel Dekker, New York, Basel, 1987.

\bibitem{Min}
N. Minorski,
Nonlinear Oscillations, Van Nostrand, New York, 1962.

\bibitem{M}
A.\,D. Myshkis,
Linear Differential Equations with Retarded Argument,
Nauka, Moscow, 1972 (in Russian).

\bibitem{N}
S.\,B. Norkin,
Differential Equations of the Second Order with Retarded Argument,
Translation of Mathematical Monographs, AMS, V. 31, Providence,
R.I., 1972.

\bibitem{Z}
B. Zhang,
On the retarded Li\'{e}nard equation,
{\em Proc. Amer. Math. Soc.} {\bf 115} (1992), no.3, 779--785.

\end{thebibliography}
\end{document}